\documentclass[reqno,a4paper,12pt]{amsart}
\usepackage[top=32mm, bottom=32mm, left=32mm, right=32mm]{geometry}
\usepackage{mathrsfs,enumerate}
\usepackage{graphicx}
\usepackage{epstopdf}
\usepackage{amssymb}
\usepackage[colorlinks=true]{hyperref}
\newtheorem{theorem}{Theorem}

\renewcommand{\descriptionlabel}[1]%
{\hspace{\labelsep}\textsf{#1}}

\theoremstyle{definition}
\newtheorem{definition}[theorem]{Definition}

\newtheorem{example}[theorem]{Example}

\newtheorem{question}[theorem]{Open Problem}

\theoremstyle{remark}
\newtheorem{remark}[theorem]{Remark}

\begin{document}

\title[Answering two open problems on Banks theorem for NADS]{Answering two open problems on Banks theorem for non-autonomous dynamical systems}


\author[X. Wu]{Xinxing Wu}
\address[X. Wu]{(Corresponding author) School of Sciences, Southwest Petroleum University, Chengdu, Sichuan, 610500, P.R. China}
\email{wuxinxing5201314@163.com}

\author[G. Chen]{Guanrong Chen}
\address[G. Chen]{Department of Electronic Engineering, City University of Hong Kong,
Hong Kong SAR, P.R. China}
\email{gchen@ee.cityu.edu.hk}

\thanks{This work was supported by the National Natural Science Foundation of China (No. 11601449),
the National Nature Science Foundation of China (Key Program) (No. 51534006), the Science and Technology
Innovation Team of Education Department of Sichuan for Dynamical Systems and its Applications (No. 18TD0013), and
the Youth Science and Technology Innovation Team of Southwest Petroleum University for Nonlinear Systems (No. 2017CXTD02).}


\date{\today}



\begin{abstract}
This paper shows that (1) there exists a topologically transitive NADS having two disjoint invariant periodic orbits with dense periodic points, which
is finitely generated but not periodic; (2) there exists a topologically transitive non-finitely generated NADS having two disjoint invariant periodic
orbits with dense periodic point, which is not sensitive. This answers positively the Open Problems 4.1 and 4.2 posed in \cite[C. Yang, Z. Li, J. Difference Equ.
Appl. 24 (2018), 1777--1782]{YL2018}.
\end{abstract}

\maketitle

Throughout, let $\mathbb{N}=\{1, 2, 3, \ldots\}$, $\mathbb{Z}^{+}=\{0, 1, 2, \ldots\}$, $\mathbb{Z}=\{\ldots, -1, 0, 1, \ldots\}$, 
and $f_{1, \infty}:=\{f_{n}\}_{n=1}^{\infty}$ be a sequence of continuous self-maps on a metric space $(X, d)$.
For any $i, n\in \mathbb{N}$, set $f_{i}^{(n)}=f_{i+(n-1)}\circ \cdots \circ f_{i}$ and $f_{i}^{0}=\mathrm{id}_{X}$.
Call $f_{1, \infty}$ a {\it non-autonomous dynamical system} (abbrev. NADS) on $X$ and denote the system by
$(X, f_{1, \infty})$. The {\it orbit} $\mathrm{orb}(x, f_{1, \infty})$ of a point $x\in X$ is the set $\{f_{1}^{(n)}(x): n\in \mathbb{Z}^+\}$, 
which is the solution of the non-autonomous difference equation
$$
\begin{cases}
x_{n+1}=f_{n+1}(x_{n}), \\
x_{0}=x.
\end{cases}
$$
In particular, if $f_{n}=f$ for all $n\in \mathbb{N}$, then the pair $(X, f)$ is a `classical' autonomous
dynamical system (abbrev. ADS).

\begin{definition}\cite{KS1996,SR2018,WZ2013,YL2018}
A NADS $(X, f_{1, \infty})$ is
\begin{enumerate}[(1)]
\item {\it topologically transitive}
for each pair non-empty open subsets $U, V$ of $X$, there exists $n\in \mathbb{Z}^+$
such that $f_1^{(n)}(U)\cap V\neq \emptyset$;
\item {\it topologically mixing} for each pair non-empty open subsets $U, V$ of $X$, there exists $N\in \mathbb{Z}^+$
such that $f_1^{(n)}(U)\cap V\neq \emptyset$ for all $n\geq N$;
\item {\it sensitive} there exists $\varepsilon>0$ such that, for any $x\in X$ and any $\delta>0$,
there exist $y\in X$ with $d(x, y)<\delta$ and $n\in \mathbb{Z}^+$ satisfying $d(f_1^{(n)}(x), f_1^{(n)}(y))\geq \varepsilon$;
\item {\it finitely generated} there exists a finite set $F$ of continuous self-maps on $X$ such that every $f_i$ of
$f_{1, \infty}$ belongs to $F$, i.e., $\{f_i: i\in \mathbb{N}\}$ is a finite set.
\end{enumerate}
\end{definition}
A point $x \in X$ is a {\it transitive point} of $f_{1, \infty}$ if its orbit $\mathrm{orb}(x, f_{1, \infty})$ is dense
in $X$. The set of all transitive points of $f_{1, \infty}$ is denoted by $\mathrm{Trans}(f_{1, \infty})$. When $\mathrm{Trans}(f_{1, \infty})=X$,
$(X, f_{1, \infty})$ is {\it minimal}. A point $x \in X$ is called an {\it equicontinuous point} if for any $\varepsilon>0$,
there exists $\delta>0$ such that whenever $y\in X$ satisfies $d(x, y)<\delta$ one has $d(f_1^{(n)}(x), f_1^{(n)}(y))<\varepsilon$ for all
$n\in \mathbb{Z}^+$. The set of all equicontinuous points of $f_{1, \infty}$ is denoted by $\mathrm{Eq}(f_{1, \infty})$. A topologically
transitive NADS is called {\it almost equicontinuous} if there exists at least one equicontinuous point. The classical Auslander-Yorke
dichotomy theorem states that every topologically transitive compact ADS is sensitive or almost equicontinuous. Recently,
we \cite[Theorem~7]{WMZL2018} extended Auslander-Yorke dichotomy theorem to Hausdorff uniform spaces.

\begin{definition}\cite{SR2018}
A point $x\in X$ is called a {\it periodic point} of $(X, f_{1, \infty})$ if there exists $n\in \mathbb{N}$
such that $f_{1}^{(nk)}(x)=x$ for all $k\in \mathbb{Z}^+$. The set of all periodic points of $(X, f_{1, \infty})$ is
denoted by $\mathrm{Per}(f_{1, \infty})$.
\end{definition}

\begin{definition}\cite{MP2013,YL2018}
Let $(X, f_{1, \infty})$ be a NADS, $x\in X$, and $A\subset X$.
\begin{enumerate}[(1)]
\item $A$ is an {\it invariant set} for $f_{1, \infty}$ if $f_{i}(A)\subset A$ for every $f_i\in f_{1, \infty}$.
\item $x$ is an {\it invariant periodic point} of $f_{1, \infty}$ if $x\in \mathrm{Per}(f_{1, \infty})$ and $\mathrm{orb}(x, f_{1, \infty})$
is an invariant set under $f_{1, \infty}$.
\end{enumerate}
\end{definition}


Banks et al. \cite{Banks1992} proved that every topologically transitive ADS with dense periodic points
is sensitive. Recently, Yang and Li \cite{YL2018} obtained the following result for generalizing Banks Theorem
from ADS to NADS, and posed two Open Problems.

\begin{theorem}\cite[Theorem~3.1]{YL2018}\label{Thm-sen}
{\em Let $X$ be a metric space without isolate points. Suppose that a finitely generated $\mathrm{NADS}$ 
$(X, f_{1, \infty})$ satisfies the following conditions:
\begin{enumerate}[(1)]
\item[(1)] $f_{1, \infty}$ is topologically transitive;
\item[(2)] $\overline{\mathrm{Per}(f_{1, \infty})}=X$;
\item[(3)] there exist two invariant periodic points $x, y\in X$ such that
$\mathrm{orb}(x,f_{1, \infty})\cap \mathrm{orb}(y, f_{1, \infty})=\emptyset$.
\end{enumerate}
Then, $(X, f_{1, \infty})$ is sensitive.}
\end{theorem}

\begin{question}\cite[Open Problem~4.1]{YL2018}\label{Q-1}
Is there any counterexample of a NADS, which is not finitely generated
but satisfies other conditions of Theorem \ref{Thm-sen}, fails to be sensitive?
\end{question}

\begin{question}\cite[Open Problem~4.2]{YL2018}\label{Q-2}
Is there any example of a NADS, which is only finitely generated but
not periodic, satisfies all the conditions of Theorem \ref{Thm-sen}?
\end{question}

This paper constructs two NADSs (Examples \ref{Exa-1} and \ref{Exa-2}) to answer positively
the above Open Problems \ref{Q-1} and \ref{Q-2}.

\begin{example}\label{Exa-1}
Fix an almost equicontinuous and non-minimal homeomorphism $(X, f)$ constructed in \cite[Theorem~4.2]{AAB1996}.
For any $n\in \mathbb{N}$, take $f_{4n-3}=f_{4n}=f^{n}$ and $f_{4n-2}=f_{4n-1}=f^{-n}$, i.e., $\{f_{n}\}_{n=1}^{+\infty}
=\{f, f^{-1}, f^{-1}, f, f^2, f^{-2}, f^{-2}, f^{2}, \ldots\}$ and let $f_{1, \infty}=\{f_{n}\}_{n=1}^{+\infty}$. Clearly, $(X, f_{1, \infty})$
is a NADS, which is not finitely generated. As $f$ is topologically transitive, it is easy to show that
$f_{1,\infty}$ is topologically transitive. For any $x\in X$, it can be verified that $\mathrm{orb}(x, f_{1, \infty})=\{f^{n}(x): n\in \mathbb{Z}\}$
and $f_{1}^{(2n)}(x)=x$ for all $n\in \mathbb{N}$. This implies that $x$ is an invariant periodic point
and $\overline{\mathrm{Per}(f_{1, \infty})}=X$.

Fix any $z\in \mathrm{Trans}(f)$. It can be shown that there exists $y\in X$ such that 
$\mathrm{orb}(z, f_{1, \infty})\cap \mathrm{orb}(y, f_{1, \infty})=\emptyset$.
In fact, suppose that for any $y\in X$, one has $\mathrm{orb}(z, f_{1, \infty})\cap \mathrm{orb}(y, f_{1, \infty})=\{f^n(z): n\in \mathbb{Z}\}\cap \{f^{n}(y): n\in \mathbb{Z}\}\neq\emptyset$. Then, there exist $n_1, n_2\in \mathbb{Z}$ such that $f^{n_1}(z)=f^{n_2}(y)$, i.e., $z=f^{n_2-n_1}(y)$, implying that
$\overline{\{f^{n}(y): n\in \mathbb{Z}^{+}\}}=\overline{\{f^{n}(z): n\in \mathbb{Z}^+\}}=X$. This means that $(X, f)$ is a minimal system, 
which is a contradiction. Thus, condition (3) in Theorem \ref{Thm-sen} is satisfied.

Applying \cite[Corollary~3.7]{AAB1996} yields that $(X, f^{-1})$ is also almost equicontinuous and
$\mathrm{Eq}(f)=\mathrm{Trans}(f)=\mathrm{Trans}(f^{-1})=\mathrm{Eq}(f^{-1})$. This implies that
$(X, f_{1, \infty})$ is not sensitive.
\end{example}

\begin{remark}
Example \ref{Exa-1} shows that there exists a NADS, which is not finitely generated
but satisfies other conditions of Theorem \ref{Thm-sen}, fails to be sensitive, giving a positive answer
to Open Problem \ref{Q-1}.
\end{remark}

\begin{example}\label{Exa-2}
Consider the space $\Sigma=\{0, 1\}^{\mathbb{N}}=\{x_1x_2 x_3\cdots: x_{i}\in \{0, 1\}, \text{ for any } i\in \mathbb{N}\}$,
with the product metric
$$
d(x, y)=\sum_{n=1}^{+\infty}\frac{|x_n-y_n|}{2^n},
$$
for any $x=x_1x_2\cdots$, $y=y_1y_2\cdots\in \Sigma$. Given any $n\in \mathbb{N}$, call $w=w_1w_2\dots w_{n}\in\{0, 1\}^n$ a
\textit{word of length $n$}, and write $|w|=n$. For any two words $u=u_1u_2\cdots u_n$ and $v=v_1v_2\cdots v_m$, call
$uv=u_1u_2\cdots u_n v_1v_2\cdots v_m$ the \textit{concatenation} of $u$ and $v$. In the same manner, call $u^m$ the
concatenation of $m$ copies of $u$ for an $m\in \mathbb{N}$ and $u^\infty$ the infinite concatenation of $u$. Let $A=a_1a_2\cdots a_{n}$
be a word. Denote $\overline{A}=\overline{a_1}\overline{a_2}\cdots \overline{a_{n}}$, and call it the {\it inverse} of $A$, where
$$
\overline{a_i}=
\begin{cases}
0, &  a_i=1, \\
1, &  a_i=0.
\end{cases}
$$
Define $\sigma: (\Sigma, d)\to (\Sigma, d)$ by
$$
\sigma(x_0 x_1 x_2 \cdots)=x_1 x_2 x_3 \cdots.
$$
Let $f_0=\sigma$ and $f_1=\sigma^{2}$. Take $A_0=0\in \{0, 1\}$, $A_1=\overline{A_0}=1\in \{0, 1\}$,
and $A_{n}=\overline{A_{0}A_1\cdots A_{n-1}}$ for all $n\geq 2$ and pick $\xi=\xi_1\xi_2\xi_3\cdots=
A_0A_1A_2\cdots A_n\cdots$, which is the classical {\it Thue-Morse sequence}. For any $i\in \mathbb{N}$,
let $g_{i}=f_{\xi_{i}}$ and $g_{1, \infty}=\{g_i\}_{i=1}^{+\infty}$.
It can be verified that $(\Sigma, g_{1, \infty})$ is
finitely generated but not periodic.

\medskip

\textsl{Claim 1.} $(\Sigma, g_{1, \infty})$ is topologically mixing.

\medskip

For any nonempty open subsets $U, V$ of $\Sigma$, since $\sigma$ is topologically mixing,
it follows that there exists $N\in \mathbb{N}$ such that for any $n\geq N$, $\sigma^{n}(U)\cap V\neq \emptyset$.
This implies that, for any $n\geq N$, $g_{1}^{(n)}(U)\cap V=f_{\xi_n}\circ \cdots \circ f_{\xi_1}(U)\cap V\neq \emptyset$; 
namely, $g_{1, \infty}$ is topologically mixing.

\medskip

\textsl{Claim 2.} $\overline{\mathrm{Per}(g_{1, \infty})}=\Sigma$.

\medskip

From the construction of $\xi$, it is clear that
$\xi_1\cdots \xi_{2|A_n|}=\overline{A_n}A_n\in \{A_{n}\overline{A_n}, \overline{A_n}A_n\}$.
$\xi_1\cdots \xi_{|A_{n+1}|}=\overline{A_n}A_n A_n \overline{A_n}$ implies that
$\xi_{2|A_n|+1}\cdots \xi_{4|A_n|}\in \{A_{n}\overline{A_n}, \overline{A_n}A_n\}$.
Then, it can be verified that for any $n\in \mathbb{N}$,
$\xi_{2j|A_{n}|+1}\cdots \xi_{2(j+1)|A_{n}|}\in \{A_{n}\overline{A_n}, \overline{A_n}A_n\}$
for all $j\in \mathbb{Z}^+$, by analogy.

For any $x=x_1x_2\cdots\in \Sigma$ and any $\delta>0$, take some $n\in \mathbb{N}$ satisfying
$\frac{1}{2^{n}}<\delta$ and let $\eta=A^{\infty}$, where $A=x_1\cdots x_{3|A_n|}$. Clearly,
$$
d(x, \eta)\leq \sum_{i=3|A_n|+1}^{+\infty}\frac{1}{2^i}=\frac{1}{2^{3|A_n|}}<\delta.
$$
From $\{\xi_{2j|A_{n}|+1}\cdots \xi_{2(j+1)|A_{n}|}: j\in \mathbb{Z}^+\}\subset \{A_{n}\overline{A_n}, \overline{A_n}A_n\}$,
it follows that for any $k\in \mathbb{Z}^+$, 
one has $g_{1}^{(2k|A_n|)}(\eta)=f_{2k|A_n|}\circ \cdots \circ f_{\xi_1}(\eta)=\sigma^{3k|A_n|}(\eta)=\eta$.
Therefore, $\eta$ is a periodic point of $g_{1, \infty}$.

\medskip

\textsl{Claim 3.} There exist two invariant periodic points, $x, y\in \Sigma$, such that
$\mathrm{orb}(x, g_{1, \infty})\cap \mathrm{orb}(y, g_{1, \infty})=\emptyset$.

\medskip

Take $x=000\cdots$ and $y=111\cdots\in \Sigma$. Clearly, $x$ and $y$ are fixed points
of $g_{1, \infty}$ and $\mathrm{orb}(x, g_{1, \infty})\cap \mathrm{orb}(y, g_{1, \infty})=\emptyset$.

\medskip

Summing up Claim 1--Claim 3 leads to the conclusion that $(\Sigma, g_{1, \infty})$ is finitely generated but
not periodic, and satisfies all the conditions of Theorem \ref{Thm-sen}. This gives a positive answer
to Open Problem \ref{Q-2}.
\end{example}


\end{document}